\documentclass[10pt]{easychair}
\usepackage{lmodern}
\usepackage[utf8]{inputenc}
\usepackage{latexsym}
\usepackage[T1]{fontenc}
\usepackage{amssymb}
\usepackage{fix-cm} 
\usepackage{stmaryrd}
\usepackage{empheq}
\usepackage{fancyhdr}
\usepackage{geometry}
\usepackage{amsmath}
\usepackage{footmisc}
\usepackage{eulervm}

\usepackage{graphicx}
\usepackage{xcolor}
\usepackage{tikz}
\usepackage{tabularx}
\usepackage{array}
\usepackage{eqparbox}
\usepackage{booktabs}

\usepackage{url}
\usepackage{hyperref}

\newcommand{\pc}{\mathbf{P}}

\title{Mereological Emptiness for the Signed Number Problem}
\author{\textbf{A. Mani}}
\institute{Indian Statistical Institute\\
203, B. T. Road, Kolkata-700108\\
\email{$a.mani.cms@gmail.com$, $amani.rough@isical.ac.in$}\\
Homepage: \url{https://www.logicamani.in}\\
Orcid: \url{https://orcid.org/0000-0002-0880-1035}}
\authorrunning{A Mani}

\titlerunning{Mereology}

\begin{document}

\maketitle


In the mathematics education literature, it is well-known that primary school students face considerable difficulties in making sense of negative integers and numbers \cite{mjg1974,kawo2010,rkkss2017,jmjf2016}. The issue was faced by mathematicians (especially western) in different perspectives for many centuries till $1850$ approximately \cite{lh2007}.  Mathematicians in China, for example, have however used signed numbers in a meaningful way \cite{wws2000,jcm1997} since ancient times. Thus, the meanings of negative numbers depend on culture and perspective, and it makes sense to construct amenable ontological explanations that are in tune with mathematically accepted meanings and use of such numbers. The purpose of this research is to address the latter problem, and propose a new way of interpreting and teaching the \emph{concepts of negative numbers} based on a mereological understanding of emptiness that has been of studied by philosophers since ancient times. In the extended version of this research, general rough sets are used to model the relation between concepts as partial algebraic systems. The focus on operations and predicates is relevant for keeping track of student reasoning and clarity. 

From the perspective of pedagogic content knowledge \cite{bdtmp2008}), \emph{approximations of real numbers} are taught in primary schools (especially in India and western countries). Relative to logicians and algebraists, students typically learn the real numbers (or rather the rationals as $\pi$ is taught as an approximation) over a permissive implicit model with a very long signature. No distinction is made between operation symbols and their interpretations (as operations), the symbol $+$ is interpreted both as a binary and a unary operation. The same is the case with the symbol $-$, while $a\div b$ may also be written as $\frac{a}{b}$. Omitting exponentiation, distinguished elements, and nullary operation symbols, the signature is $\Sigma = (\leq, \geq, +, \times , \div,  -, /, \oplus, \ominus, (2, 2, 2, 2, 2, 2, 2, 1, 1) )$. The numbers in braces refer to the place value of operation symbols, and the \emph{plus and minus signs} have been denoted as $\oplus$ and $\ominus$ respectively. \emph{The standard conception of terms in algebra (in higher classes) is partially split into a number of concepts such as multiplicative terms, binomials, monomials, rational fractions, and fractions in school algebra}. Multiplicative terms are also simply referred to as terms, and even these are not provided with a recursive definition.   

The predicate $\leq$ is used in the sense of \emph{being numerically less than} is a transitive parthood. In the context of word problems, ideas of some objects being bigger or smaller than others often leads to ambiguity and vagueness especially when the objects have multiple attributes associated. Functional non-transitive parthood is apparently not commonly used in mathematics text books at least -- though this often happens in the context of activities and descriptions that are very contextual and under-specified. So the addition of an additional parthood predicate $\pc$ to the signature is a good idea. The partial algebraic systems that are taught to school children are therefore related to $\mathfrak{R} = \left\langle\underline{A}^*, \,\Sigma  \right\rangle$ with $\underline{A}^*$ being the set corresponding to the algebraic closure of the union of algebraic numbers and $\{\pi\}$. The axioms/properties satisfied are not uniformly learnt relative to concept formation/maturation levels. In addition, \emph{alternate conceptual axioms} (possibly higher order) may be assumed and all this can be used to explain only aspects of diverse problem-solving behaviors, and concepts of negative numbers adhered to. 

Students often fail to understand the negative integers, fractions and numbers because they cannot effectively associate them with anything material in the real world. It is also that most people understand negative real numbers and elements of a vector space (or linear algebra) relative  to ideas of graphical representation, reference frames, states, and transformation or change.  For example, in \cite{kawo2010}, the authors suggest that the signed integer $-2$ can be interpreted \emph{as a change}, \emph{as a state} or \emph{as a directed relation between numbers and quantities}. Visualizations based on ideas of displacement \cite{rkkss2017} relate to a directed relation between numbers and quantities. In the mentioned paper, the authors refer to displacement of a lift in upward and downward directions relative to ground level (level zero) as a means of introducing the idea of negative numbers. Such ideas cannot be naturally transported to other material contexts, and are suited for handling small negative integers alone. While models based on neutralization (like equal $+v$ and $-ve$ charges may cancel each other to produce a net $0$ charge) lead to ontological inconsistencies that can hinder learning.

For an introduction to mereology, the reader is referred to \cite{ham2017,rgac15,am2021}. Here non-transitive mereology will be adhered to with scope for functional versions of parthood and vagueness \cite{am5586,am9969,seibtj2015}. Aspects of mereology are always used in the mathematical discourse. Because the formalisms of school mathematics are not uniformly consistent, only fragments can be assumed to be well formulated. Language use in the mathematics class can involve both inconsistency and vagueness. But whatever ideas of \emph{bigger than}, \emph{smaller than}, \emph{operations being connected with order or parthood}, \emph{composed of}, \emph{being an ingredient of}, \emph{being connected} and other mereological predicates  are part of the discourse in \emph{some form}. A student of primary school may not have a proper conception of connected parts especially from a geometrical perspective, or some other type. 

In the present author's view, the problem of signs is due to inconsistencies in the ascription of ontologies to the real numbers. Specifically the concept of zero is read in a way compatible with a field (in the sense of algebra) alone. Most mereo-ontologies are not compatible with this idea of zero associated with a universal concept of emptiness. In fact almost every mereological theory may be associated with ideas of typed emptiness. This holds for both ancient and more recent formal mereologies. If the \emph{content of perception} is seen as distinct from \emph{the result of inference through a logical process}, then entirely different concepts of emptiness are possible - an emptiness would be a perceived object in the former, while the latter would possibly yield a relatively absolute nothingness. Examples of the former can be found even in Huayan texts, where ontological rules are specified in addition.

\emph{Most importantly, it suffices to introduce a large number of common possible types of emptiness to enable concept development}. Thus what is traditionally taught as \emph{an empty glass} can be a $-50$ \emph{marbles empty glass} if the glass in question can be filled with $50$ marbles and it is  that the proposed activity concerns filling glasses with marbles. At the same time, the glass is not $-1$ \emph{cabbage empty} because it cannot fit in a cabbage, though it is $-1$ drink empty. The idea of zero in these contexts would then be about possible give and take operations in the absence of such in reality. It may be noted that in ancient Chinese mathematics, a debt would be the same as possessing a negative sum.

The proposed approach would also help in highlighting the difference between numbers used for counting, and the thing being counted. A single pencil \emph{has the property of being one}, two pencils \emph{that of being two}, and  so on. This meaning differs from the intention of a mark of \emph{one} or \emph{two} on a ruler.  Related to all this is the question \emph{Why should the oneness of a pencil be the same as that of an eraser}? The principle of equivalence associated is not absolutely essential in any actual mathematical application or model. If the solution of a partial differential equation modeling the flow of a fluid is $x=1, \, y= 1$ and $z=x$, then does it ever mean that anything with oneness can be freely substituted into the interpretation or meaning of the solution?

At the same time, the apparent ambiguity and vagueness induced by the approach relative to universal ideas of emptiness, suggest that temporal aspects should be integrated into ideas such as \emph{a glass is $-1$ drink empty}. $-$ in $-1$ here is a reference to the ability of the glass to potentially contain a drink in future. The concept of \emph{a drink} is intrinsically vague from a measure-theoretic point of view -- that is besides the scope of this discourse. The word \emph{empty} can then be understood as a universal vague reference to the property of having some context-dependent $-$ sign. The alternate interpretation of the empty and of zeros naturally raises questions about operations involving the latter. More details may be found in the extended version of this research  by the present author.

\begin{small}
\bibliographystyle{plain.bst}
\bibliography{algroughf69fl.bib}
\end{small}
\end{document}